\documentclass[letterpaper,10pt]{amsart}
\usepackage{amsmath,amssymb,amsxtra,amsthm}
\usepackage{graphicx}
\usepackage[usenames]{color}

\theoremstyle{plain}
\newtheorem{theorem}{Theorem}[section]

\newtheorem{lemma}[theorem]{Lemma}

\theoremstyle{definition}

\newtheorem{remark}[theorem]{Remark}

\newcommand{\B}{\mathbb}

\newcommand{\z}{\zeta}
\newcommand{\ga}{\alpha}

\newcommand{\eps}{\varepsilon}

\newcommand{\gl}{\lambda}

\newcommand{\gq}{\vartheta}

\newcommand{\gz}{\zeta}

\begin{document}
\title[Residue class distribution of $\Omega(n)$]{On the residue class distribution of the number of prime divisors of an integer}
\author{Michael Coons}
\author{Sander R. Dahmen}
\address{Simon Fraser University, Burnaby, British Columbia, Canada, V5A 1S6}
\email{mcoons@sfu.ca, sdahmen@irmacs.sfu.ca}
\date{\today}
%%%%%%%%%%%%%%%%%%%%%%%%%%%%%%%%%%%%%%%%%%%%%%%%%%%%%%%%%%%%%%%%%%%

\begin{abstract} The {\em Liouville function} is defined by $\gl(n):=(-1)^{\Omega(n)}$ where $\Omega(n)$ is the number of prime divisors of $n$ counting multiplicity. Let $\z_m:=e^{2\pi i/m}$ be a primitive $m$--th root of unity. As a generalization of Liouville's function, we study the functions $\gl_{m,k}(n):=\z_m^{k\Omega(n)}$. Using properties of these functions, we give a weak equidistribution result for $\Omega(n)$ among residue classes. More formally, we show that for any positive integer $m$, there exists an $A>0$ such that for all $j=0,1,\ldots,m-1,$ we have  $$\#\{n\leq x:\Omega(n)\equiv j\ (\bmod\ m)\}=\frac{x}{m}+O\left(\frac{x}{\log^A x}\right).$$ Best possible error terms are also discussed. In particular, we show that for $m>2$ the error term is not $o(x^\ga)$ for any $\ga<1$.
\end{abstract}

\maketitle

%%%%%%%%%%%%%%%%%%%%%%%%%%%%%%%%%%%%%%%%%%%%%%%%%%%%%%%%%%%%%%%%%%%
\section{Introduction}
%%%%%%%%%%%%%%%%%%%%%%%%%%%%%%%%%%%%%%%%%%%%%%%%%%%%%%%%%%%%%%%%%%%

The {\em Liouville function}, denoted $\gl(n)$, is defined by $\gl(n):=(-1)^{\Omega(n)}$ where $\Omega(n)$ is the number of prime divisors of $n$ counting multiplicity. The Liouville function is intimately connected to the Riemann zeta function and hence to many results and conjectures in prime number theory. Recall that \cite[pp.~617--621]{Lan1} for $\Re s>1$, we have \begin{equation}\label{lquo} \sum_{n\geq 1}\frac{\gl(n)}{n^s}=\frac{\gz(2s)}{\gz(s)},\end{equation} so that $\gz(s)\neq 0$ for $\Re s\geq \gq$ provided $$\sum_{n\leq x}\gl(n)=o(x^\gq).$$ The prime number theorem allows the value $\gq=1$, so that for $j=0,1$ we have that $$\#\{n\leq x:\Omega(n)\equiv j\ (\bmod\ 2)\}\sim \frac{x}{2}.$$ We generalize this result to the following theorem.

\begin{theorem}\label{main} Let $m$ be a positive integer and $j=0,1,\ldots,m-1$. Then the (natural) density of the set of all $n\in\B{Z}_{> 0}$ such that $\Omega(n)\equiv j\ (\bmod\ m)$ exists, and is equal to $1/m$; furthermore, there exists an $A>0$ such that  $$N_{m,j}(x):=\#\{n\leq x:\Omega(n)\equiv j\ (\bmod\ m)\}=\frac{x}{m}+O\left(\frac{x}{\log^A x}\right).$$
\end{theorem}

In order to prove this theorem, we study a generalization of Liouville's function. Namely, let $m$ be a positive integer and $\z_m:=e^{2\pi i/m}$ be a primitive $m$--th root of unity. Define $$\lambda_{m,k}(n):=\z_m^{k\Omega(n)}.$$ As with $\gl(n)$, since $\Omega(n)$ is completely additive, $\gl_{m,k}(n)$ is completely multiplicative. For $\Re s>1$, denote $$L_{m,k}(s):=\sum_{n\geq 1}\frac{\gl_{m,k}(n)}{n^{s}}.$$

The functions $\gl_{m,k}(n)$ and $L_{m,k}(s)$ were introduced by Kubota and Yoshida \cite{Kub1}. They gave (basically) a multy-valued analytic continuation of $L_{m,k}(s)$ to the region $\Re s>1/2.$ Using this, for $m\geq 3$ and $k=1,\ldots,m-1$ with $m/k\neq 2$, they showed that certain asymptotic bounds on the partial sums $$S_{m,k}(x):=\sum_{n\leq x}\gl_{m,k}(n),$$ cannot hold; in particular, this sum cannot be $o(x^\ga)$ for any $\ga<1$. 
Finally, this is used by the authors to show, given Theorem \ref{main}, that if $m\geq 3$, then an asymptotic of the form 
\begin{equation}\label{asymptotic}
N_{m,j}(x)=\frac{x}{m}+o(x^\ga)
\end{equation}
cannot hold simultaneously for all $j=0,1,\ldots,m-1,$ if $\ga<1$. We will show that if $m \geq 3$, then for \emph{all} $j=0,1,\ldots,m-1$ the asymptotic \eqref{asymptotic} does not hold if $\alpha<1$.
This is in striking contrast to the expected result for $m=2$. Recall that in the case that $m=2$, if the Riemann hypothesis is true then $$ N_{2,j}(x)=\frac{x}{2}+o(x^{1/2+\eps})$$ for $j=0,1$ and any $\eps>0$.

%%%%%%%%%%%%%%%%%%%%%%%%%%%%%%%%%%%%%%%%%%%%%%%%%%%%%%%%%%%%%%%%%%%
\section{Perliminary results}
%%%%%%%%%%%%%%%%%%%%%%%%%%%%%%%%%%%%%%%%%%%%%%%%%%%%%%%%%%%%%%%%%%%

\begin{lemma}\label{qOmega} Let $m$ be a positive integer. Then
for $k=0, 1, \ldots, m-1,$ we have
\begin{equation}\label{S}
S_{m,k}(x) = \sum_{j=0}^{m-1} \z_m^{jk} N_{m,j}(x)
\end{equation}
and for $j=0, 1, \ldots, m-1,$ we have
\begin{equation}\label{N}
N_{m,j}(x) = \frac{1}{m}\sum_{k=0}^{m-1}\z_m^{-jk}S_{m,k}(x).
\end{equation}
\end{lemma}

\begin{proof}
We have
\begin{align*}
S_{m,k}(x) & = \sum_{n \leq x} \z_m^{k \Omega(n)}\\
 & =  \sum_{j=0}^{m-1} \sum_{\substack{n \leq x\\ \Omega(n) \equiv j\ ({\rm mod}\ m)}}\z_m^{k \Omega(n)}\\
 & =   \sum_{j=0}^{m-1} \z_m^{kj} N_{m,j}(x),
\end{align*}
which proves the first formula of the lemma. Instead of directly inverting the matrix determined by this formula, we proceed as follows to obtain the second formula.
Consider the right--hand side of \eqref{N}. Using the definition of $\gl_{m,k}(n)$ we have $$\frac{1}{m}\sum_{k=0}^{m-1}\z_m^{-jk}\sum_{n\leq x}\gl_{m,k}(n)=\sum_{n\leq x}\frac{1}{m}\sum_{k=0}^{m-1} \z_m^{(\Omega(n)-j)k}.$$ If $n$ satisfies $\Omega(n)\equiv j\ (\bmod\ m)$, then $\z_m^{\Omega(n)-j}=1$, so that $$\frac{1}{m}\sum_{k=0}^{m-1} \z_m^{(\Omega(n)-j)k}=1.$$ If $n$ does not satisfy $\Omega(n)\equiv j\ (\bmod\ m)$, then $\z_m^{\Omega(n)-j}\neq 1$. We thus have $$\frac{1}{m}\sum_{k=0}^{m-1} \z_m^{(\Omega(n)-j)k}=\frac{1}{m}\cdot\frac{\z_m^{(\Omega(n)-j)m}-1}{\z_m^{(\Omega(n)-j)}-1}=\frac{1}{m}\cdot\frac{0}{\z_m^{(\Omega(n)-j)}-1}=0.$$ This proves the second part of the lemma.
\end{proof}

To yield our density result on the number of prime factors, counting multiplicity, modulo $m$, we need the following result.

\begin{theorem}\label{Tglqk0} For every $m\in\B{Z}_{> 0}$ there is an $A>0$ such that for all $k=1,$ $\ldots,m-1,$ we have $$\left| S_{m,k}(x) \right|\ll \frac{x}{\log^{A} x}.$$
\end{theorem}

\noindent To prove this, we use the following theorem.

\begin{theorem}[Hall \cite{Hall}]\label{hall} Let $D$ be a convex subset of the closed unit disk in $\B{C}$ containing $0$ with perimeter $L(D)$. If $f:\B{Z}_{> 0}\to\B{C}$ is a multiplicative function with $|f(n)|\leq 1$ for all $n\in\B{Z}_{>0}$ and $f(p)\in D$ for all primes $p$, then \begin{equation}\label{1.5} \frac{1}{x}\left|\sum_{n\leq x}f(n)\right|\ll \exp\left(-\frac{1}{2}\left(1-\frac{L(D)}{2\pi}\right)\sum_{p\leq x}\frac{1-\Re f(p)}{p}\right).\end{equation}
\end{theorem}

\begin{proof} This is a direct consequence of Theorem 1 of \cite{Hall}.
\end{proof}

\begin{proof}[Proof of Theorem \ref{Tglqk0}] Set $D$ equal to the convex hull of the $m$--th roots of unity and $f=\gl_{m,k}$ . Because $D$ is a convex subset strictly contained in the closed unit disk of $\B{C}$, we have $L(D)<2\pi$. This gives $$c:=\frac{1}{2}\left(1-\frac{L(D)}{2\pi}\right)>0.$$ Applying Theorem \ref{hall} yields \begin{align*} \frac{1}{x}\left|\sum_{n\leq x}\gl_{m,k}(n)\right| &\ll \exp\left(-c\sum_{p\leq x}\frac{1-\Re \gl_{m,k}(p)}{p}\right)\\
&= \exp\left(-c(1-\Re\z_m^k)\sum_{p\leq x}\frac{1}{p}\right)
\end{align*} Since $\sum_{p\leq x}p^{-1}=\log\log x+O(1)$, this quantity is \begin{align*} &\ll \exp\left(-c(1-\Re\z_m^k)\log\log x\right) \\ &= \left(\frac{1}{\log x}\right)^{c(1-\Re\z_m^k)}.\end{align*} Noting that $0<k<m$, we have $c(1-\Re\z_m^k)>0$. Set $A:=\min_{0<k<m} \{c(1-\Re\z_m^k)\}$. Then $A>0$ and we obtain \begin{equation*}\left|\sum_{n\leq x}\gl_{m,k}(n)\right|\ll \frac{x}{\log^{A} x}.\qedhere\end{equation*}
\end{proof}

%%%%%%%%%%%%%%%%%%%%%%%%%%%%%%%%%%%%%%%%%%%%%%%%%%%%%%%%%%%%%%%%%%%
\section{Proof of Theorem \ref{main}}
%%%%%%%%%%%%%%%%%%%%%%%%%%%%%%%%%%%%%%%%%%%%%%%%%%%%%%%%%%%%%%%%%%%

\begin{proof}[Proof of Theorem \ref{main}]  Lemma \ref{qOmega} directly gives us 
\begin{equation}\label{errorterm}
N_{m,j}(x)=\frac{1}{m}S_{m,0}(x)+\frac{1}{m}\sum_{k=1}^{m-1}\z_m^{-jk} S_{m,k}(x).
\end{equation}
The first term of the right-hand side \eqref{errorterm} is $$\frac{1}{m} S_{m,0}(x)=\frac{1}{m}\sum_{n\leq x} 1=\frac{x}{m}+o(1).$$ Applying the triangle inequality and Theorem \ref{Tglqk0}, we get that the absolute value of the second term of the right-hand side of \eqref{errorterm} is $$\left| \frac{1}{m}\sum_{k=1}^{m-1}\z_m^{-jk} S_{m,k}(x) \right| \leq \frac{1}{m}\sum_{k=1}^{m-1}\left| S_{m,k}(x) \right|\ll \frac{x}{\log^A x}$$ for some $A>0$.
This gives us our desired result.\end{proof}

%%%%%%%%%%%%%%%%%%%%%%%%%%%%%%%%%%%%%%%%%%%%%%%%%%%%%%%%%%%%%%%%%%%
\section{Results for error terms}
%%%%%%%%%%%%%%%%%%%%%%%%%%%%%%%%%%%%%%%%%%%%%%%%%%%%%%%%%%%%%%%%%%%

For $m \in \B{Z}_{> 0}$ and $j=0, 1, \ldots, m-1,$ we introduce the error term $$R_{m.j}(x):=N_{m,j}(x)-\frac{x}{m}.$$ Theorem \ref{main} implies that $$R_{m,j}(x)=o(x).$$ For $m>2$, Kubota and Yoshida \cite{Kub1} prove, conditionally on Theorem \ref{main}, that at least one of the error terms $R_{m,j}(x)$ is not $o(x^{\ga})$ for any $\ga<1$. We strengthen their result (unconditionally) as follows.

\begin{theorem}\label{error}
Let $m \in \B{Z}_{> 2}$ and let $\alpha<1$. Then none of $R_{m,0}, R_{m,1}, \ldots,$ $R_{m,m-1}$ are $o\left(x^{\alpha}\right)$.
\end{theorem}

Following \cite{Kub1}, we use the following results.

\begin{lemma}\label{perron}
Let $\{a_n\}_{n \in \B{Z}_{>0}}$ be a sequence of complex numbers and $\alpha>0$. If the partial sums satisfy $\sum_{n \leq x}a_n=o\left(x^{\alpha}\right),$ then the Dirichlet series $\sum_{n \geq 1} a_n n^{-s}$ converges for $\Re s>\alpha$ to a holomorphic (single--valued) function.
\end{lemma}

\begin{proof}
This follows directly from Perron's formula \cite[p.~243 Lemma 4]{Apostol}.
\end{proof}

\begin{theorem}\label{dircon}
Let $m \in \B{Z}_{> 2}$ and let $k=1, 2, \ldots, m-1$. The Dirichlet series $L_{m,k}(s)$ can be analytically continued to a multi--valued function on $\Re s>1/2$ given by the product $\gz(s)^{\z_m^k}G_{m,k}(s)$ where $G_{m,k}(s)$ is an analytic function for $\Re s>1/2.$ In particular, if $k\not=m/2$, then for any $\ga<1$ the Dirichlet series $L_{m,k}(s)$ does not converge for all $s$ with $\Re s>\ga$.
\end{theorem}

\begin{proof}
The first part follows from Theorem 1 in \cite{Kub1} (strictly speaking this handles only the case $k=1$, but the proof of this theorem works for general $k$). Note that $\z_m^k$ is not rational for $k\not=m/2$. Since $\zeta(s)$ has a pole at $s=1$, this mens that no branch of $\zeta(s)^{\z_m^k}$ is holomorphic in a neighbourhood of $s=1$.
\end{proof}

Let $m>2$ and $j=0, 1, \ldots, m-1$. From \eqref{N} we get
$$R_{m,j}(x)=\frac{1}{m}\sum_{k=1}^{m-1}\z_m^{-jk}S_{m,k}(x)-\frac{\{x\}}{m},$$ where $\{x\}$ denotes the fractional part of $x$.
In light of Lemma \ref{perron}, to obtain that $R_{m,j}(x)$ is not $o(x^{\ga})$ for any $\alpha<1$, it suffices to show that the generating function of $R_{m,j}(x)+\{x\}/m$, which is $$\sum_{k=1}^{m-1} \z_m^{-jk}
L_{m,k}(s),$$ cannot be analytically continued to a holomorphic (single--valued)
function in the half plane $\Re s > \alpha$.
\begin{remark}
We can quickly obtain the result for at least two of the error terms as follows. For $k=1, 2, \ldots, m-1$, using \eqref{S} we have $$S_{m,k}(x)=\sum_{j=0}^{m-1}\z_m^{jk} R_{m,j}(x).$$ By Lemma \ref{perron} and Theorem \ref{dircon}, $S_{m,k}(x)$ is not $o\left(x^{\alpha}\right)$ for any $\alpha<1$, so that at least one of the error terms $R_{m,j}(x)$ is not $o\left(x^{\alpha}\right)$, which is the above mentioned result of Kubota and Yoshida. From (\ref{S}) with $k=0$, we obtain $$\sum_{j=0}^{m-1}R_{m,j}(x)=S_{m,0}(x)-x=-\{x\}.$$
This shows that it is impossible that all but one of the error terms $R_{m,j}(x)$ are $o\left(x^{\alpha}\right)$ for an $\alpha<1$.
\end{remark}

We now proceed with the proof of the main result of this section.

\begin{proof}[Proof of Theorem \ref{error}]
Let $1/2<\alpha<1$ and let $c_1,c_2,\ldots,c_{m-1} \in \B{C}^*$.
We shall prove that the linear combination
\[f(s):=\sum_{k=1}^{m-1} c_k L_{m,k}(s)\]
cannot be analytically continued to a holomorphic (single--valued)
function in the half plane $\Re s > \alpha$. Suppose to the contrary that it can and
assume for now that $L_{m,1}(s),L_{m,2}(s), \ldots, L_{m,m-1}(s)$ are
linearly independent over $\B{C}$, which shall be shown later. Let $C$ denote a closed loop in
the half plane $\Re s > \alpha$ winding around $s=1$ once in the
positive direction and not around any zeroes of $\zeta(s)$. As
pointed out in \cite{Kub1}, the analytic continuation of $L_{m,k}(s)$
along $C$ gives us $\exp \left( -2 \pi i \z_m^k \right) L_{m,k}(s)$.
From the holomorphicity assumption on $f(s)$, it follows that the analytic
continuation of $f(s)$ along $C$ is $f(s)$ itself. So
$$\sum_{k=1}^{m-1} c_k L_{m,k}(s)=\sum_{k=1}^{m-1} c_k \exp \left( -2 \pi i
\z_m^k\right) L_{m,k}(s),$$ and from the linear independence over
$\B{C}$ of the functions $L_{m,k}(s)$, we obtain that $\exp \left( -2 \pi i
\z_m^k \right)=1$ for $k=1,2,\ldots, m-1$. This means $\z_m^k \in
\B{Z}$ for $k=1,2,\ldots, m-1$, a contradiction if $m>2$.

We are left with proving that $L_{m,1}(s),L_{m,2}(s), \ldots, L_{m,m-1}(s)$
are linearly independent over $\B{C}$. This can be done along
similar lines. Suppose they are not linearly independent over
$\B{C}$. Let $b$ be the smallest integer such that there exists a
nontrivial linear dependence over $\B{C}$ of $b$ different functions
$L_{m,k}(s)$, say $L_{m,k_1}(s),L_{m,k_2}(s),\ldots, L_{m,k_b}(s)$ for
$0<k_1<k_2< \ldots < k_b<m$. Since the function $L_{m,k}(s)$ are nonzero, we
have $b \geq 2$, furthermore
$$L_{m,k_1}(s)=\sum_{n=2}^b d_n L_{m,k_n}(s)$$
for some $d_2,\ldots,d_b \in \B{C}^*$. Analytic continuation along
$C$ yields
\begin{align*}
\exp \left(-2 \pi i \z_m^{k_1}\right) L_{m,k_1}(s) = & \sum_{n=2}^b
d_n
\exp \left( -2 \pi i \z_m^{k_1} \right) L_{m,k_n}(s)\\
 = & \sum_{n=2}^b d_n
\exp \left( -2 \pi i \z_m^{k_n} \right) L_{m,k_n}(s).
\end{align*}
By the minimality of $b$, we have that the $b-1 \geq 1$ functions
$L_{m,k_2}(s),\ldots, L_{m,k_b}(s)$ are linearly independent over $\B{C}$,
so $\exp \left( -2 \pi i \z_m^{k_1} \right)=\exp \left( -2 \pi i
\z_m^{k_n} \right)$ for $n=2,\ldots,b$. This means
$\z_m^{k_1}-\z_m^{k_n} \in \B{Z}$ for $n=2,\ldots,b$. One easily
obtains that the only possibility for this is when $b=2$ and
$(\z_m^{k_1},\z_m^{k_2})=(1/2 + 1/2 \sqrt{-3},-1/2 + 1/2
\sqrt{-3})$ or $(\z_m^{k_1},\z_m^{k_2})=(-1/2 - 1/2 \sqrt{-3},1/2
- 1/2 \sqrt{-3})$. Therefore, to complete the proof of the independence result, it suffices to show that $L_{6,1}(s)/L_{6,2}(s)$
and $L_{6,4}(s)/L_{6,5}(s)$ are not constant. To see this, we use the formula $L_{m,k}(s)=\gz(s)^{\z_m^k}G_{m,k}(s)$, which readily gives $$\frac{L_{6,1}(s)}{L_{6,2}(s)}=\gz(s)\frac{G_{6,1}(s)}{G_{6,2}(s)}.$$ The function $\gz(s)$ has a pole at $s=1$ and $G_{m,k}(1)\neq 0$, since for $\Re s>1/2$ we have $$\prod_{k=1}^{m-1} G_{m,k}(s)=\gz(ms).$$ We conclude that $L_{6,1}(s)/L_{6,2}(s)$ is not constant. The proof of the result for   $L_{6,4}(s)/L_{6,5}(s)$ follows similarly. This completes the proof.
\end{proof}

\begin{remark} In the spirit of prime numbers races, it seems fitting that further study should be taken to investigate the sign changes of $N_{m,j}(x)-N_{m,j'}(x)$ for $j\neq j'.$ For the case $m=2$ some such investigations have been undertaken; see \cite{BFM} and the references therein.
\end{remark}

%%%%%%%%%%%%%%%%%%%%%%%%%%%%%%%%%%%%%%%%%%%%%%%%%%%%%%%%%%%%%%
%%%%%%%%%%%%%%%%%%%%%%%%%%%%%%%%%%%%%%%%%%%%%%%%%%%%%%%%%%%%%%

\bibliographystyle{amsplain}
\providecommand{\bysame}{\leavevmode\hbox to3em{\hrulefill}\thinspace}
\providecommand{\MR}{\relax\ifhmode\unskip\space\fi MR }
\providecommand{\MRhref}[2]{%
  \href{http://www.ams.org/mathscinet-getitem?mr=#1}{#2}
}
\providecommand{\href}[2]{#2}

\end{document}